\documentclass[a4, 12pt]{amsart}

\usepackage[mathscr]{eucal}
\usepackage{amssymb}
\usepackage{latexsym}
\usepackage{amsthm}
\usepackage{amsmath,amssymb}
\usepackage{type1cm}
\usepackage{amsfonts}
\usepackage{comment}

\makeatletter
\@addtoreset{equation}{section}
\makeatother

\newtheorem{Def}{Definition}[section]
\newtheorem{Thm}{Theorem}[section]

\newtheorem{Rem}{Remark}[section]

\newtheorem{Lem}{Lemma}[section]
\newtheorem{Ex}{Example}[section]

\newcommand{\fin}{
\begin{flushright}
$ \Box $
\end{flushright}
}
\newcommand{\innprod}[2]{\left\langle #1,~ #2 \right\rangle}

\newcommand{\ep}{\varepsilon}

\newcommand{\cL}{{\mathcal L}}

\newcommand{\dm}{d\mu}

\title[Complete $\lambda$-hypersurfaces]
{A global pinching theorem  of  complete $\lambda$-hypersurfaces}
\author[S. Ogata]{Shiho Ogata}

\begin{document}
\maketitle

\begin{abstract}
In this paper,  the pinching problems of complete $\lambda$-hypersurfaces in a Euclidean space 
$\mathbb R^{n+1}$ are studied. By making  use of the Sobolev inequality, we prove a global pinching 
theorem of complete $\lambda$-hypersurfaces
in a Euclidean space  $\mathbb R^{n+1}$. 
\end{abstract}


\section{Introduction}

Let $ M^n $ be an $n $-dimensional manifold, and $X : M^n \rightarrow \mathbb{R}^{n+1} $  an immersed hypersurface in 
a Euclidean space $ \mathbb{R}^{n+1} $.  If  $ X : M^n \rightarrow \mathbb{R}^{n+1} $ satisfies 
\begin{eqnarray*}
H + \innprod{ X }{ N } = 0,
\end{eqnarray*}
one calls  that $ X  : M^n \rightarrow \mathbb{R}^{n+1} $ is a  self-shrinker of mean curvature flow,  
where $ H $ and $ N $ are the mean curvature and the unit normal vector of $ X : M^n \rightarrow \mathbb{R}^{n+1} $, 
respectively, and $ \innprod{\cdot}{\cdot} $ is the standard  inner product of  $ \mathbb{R}^{n+1}$.
\begin{Rem}
If $ X : M^n \rightarrow \mathbb{R}^{n+1} $ is a  self-shrinker of mean curvature flow, then $X(t)=\sqrt {1-2t}X$ is a self-similar solution 
 of mean curvature flow. 
\end{Rem}
It is well-known that  the $ n $-dimensional Euclidean space $ \mathbb{R}^n$, 
the $ n $-dimensional sphere $ S^n ( \sqrt{n} ) $ 
and the $ n $-dimensional cylinder $ S^k (\sqrt{k} ) \times \mathbb{R}^{n-k}$, for $ 1 \leq k \leq n-1 $, 
 are  the standard  self-shrinkers of mean curvature flow. For the other examples of self-shrinkers of  mean curvature flow, see \cite{A}, \cite{D}, 
\cite{KKM}, \cite{KM} and \cite{M}.
 
$ X (t) : M^n \rightarrow \mathbb{R}^{n+1}$ is called a variation of $ X: M^n \rightarrow \mathbb{R}^{n+1} $ if $ X (t) : M^n \rightarrow \mathbb{R}^{n+1} , t \in ( -\ep , \ep ) $,  are  a  family of immersions with $ X(0) = X $.  
We define a weighted area functional $ A : ( -\ep , \ep ) \rightarrow \mathbb{R} $ as follows:
\begin{eqnarray*}
A(t) = \int_M e^{-\frac{|X(t)|^2 }{2}} \dm_t ,
\end{eqnarray*} 
where $ \dm_t $ is the area element of $ X (t) : M^n \rightarrow \mathbb{R}^{n+1}$. 
In \cite{CM}, Colding and Minicozzi  have proved that that  $ X : M^n \rightarrow \mathbb{R}^{n+1} $ is a critical point of the weighted area functional $ A(t) $  if and only if $ X : M^n \rightarrow \mathbb{R}^{n+1} $ is a self-shrinker of mean curvature flow. 

In \cite{CL}, Cao and Li (cf. \cite{CM}, \cite{CP} \cite{LW} and \cite{S}) have proved a gap theorem of complete self-shrinkers of mean curvature flow as follows: 

\begin{Thm}
Let $ X : M^n \rightarrow \mathbb{R}^{n+1} $ be  an $n$-dimensional complete proper self-shrinker in $ \mathbb{R}^{n+1} $. 
If the squared norm $S$ of the second fundamental form of $ X : M^n \rightarrow \mathbb{R}^{n+1} $ satisfies $ S \leq 1 $, then $ X : M^n \rightarrow \mathbb{R}^{n+1} $  is isometric to one of the following: 
\begin{enumerate}
\item the sphere $ S^n(\sqrt{n}) $, 
\item the Euclidean space $\mathbb{R}^n$, 
\item the cylinder $ S^k (\sqrt{k}) \times \mathbb{R}^{n-k} $ for $1 \leq k \leq n-1$.
\end{enumerate}
\end{Thm}

By using the following  Sobolev inequality for $n$-dimensional complete hypersurfaces: 
\begin{eqnarray*}
\kappa^{-1} \left( \int_M g^\frac{2n}{n-2} \dm \right)^\frac{n-2}{n} \leq \int_M |\nabla g|^2 \dm + \frac{1}{2} \int_M H^2 g^2 \dm ,~~~ ^\forall g \in C_c^\infty (M),
\end{eqnarray*}
where $ \kappa > 0 $ is a constant, 
 Ding and Xin \cite{DX} have proved  a  rigidity theorem of complete self-shrinkers of mean curvature flow as follows: 

\begin{Thm}
Let $ X : M^n \rightarrow \mathbb{R}^{n+1} $ be an $n$-dimensional complete immersed self-shrinker  of mean curvature flow in $ \mathbb{R}^{n+1} $. 
If $ X : M^n \rightarrow \mathbb{R}^{n+1} $ satisfies 
\begin{eqnarray*}
\left( \int_M S^\frac{n}{2} \dm \right)^\frac{2}{n} < \dfrac{4}{3n\kappa},
\end{eqnarray*}
then $ X : M^n \rightarrow \mathbb{R}^{n+1} $ is isometric to the Euclidean space $ \mathbb{R}^n $, where $ S $ denotes the squared norm
of the second fundamental form of $ X : M^n \rightarrow \mathbb{R}^{n+1} $.  
\end{Thm}

In \cite{CW}, Cheng and Wei  have  introduced a notation of so-called $ \lambda $-hypersurfaces of the weighted volume-preserving mean curvature as follows: 

\begin{Def}
Let $ X : M^n \rightarrow \mathbb{R}^{n+1} $  be an $n$-dimensional immersed hypersurface in $ \mathbb{R}^{n+1} $.  If 
\begin{eqnarray*}
H + \innprod{ X}{N} = \lambda
\end{eqnarray*}
is satisfied,  
where $ \lambda $ is constant, then  $ X : M^n \rightarrow \mathbb{R}^{n+1} $ is  called a $ \lambda $-hypersurface of the weighted volume-preserving mean curvature.  For simple, we call it a $ \lambda $-hypersurface.
\end{Def}

\begin{Rem} 
From definition, we know that  if $ \lambda = 0 $, $ X : M^n \rightarrow \mathbb{R}^{n+1} $ is a self-shrinker of mean curvature flow. 
\end{Rem}

\begin{Ex}
All of self-shrinkers of mean curvature flow is  $ \lambda $-hypersurfaces with $ \lambda = 0 $. 
\end{Ex}

\begin{Ex}
The $ n $-dimensional sphere $ S^n (r ) $ with $ r > 0 $ is a compact $ \lambda $-hypersurface with 
$ \lambda =  \frac{n}{r} - r $. We should notice that only the  $n$-dimensional sphere $S^n(\sqrt n)$ is the self-shrinker of mean curvature 
flow.
\end{Ex}

\begin{Ex}
The $ n $-dimensional cylinder $ S^k (r ) \times \mathbb{R}^{n-k} $ with $ r > 0 $ for $ 1 \leq k \leq n-1 $ is a complete and non-compact $ \lambda $-hypersurface with $ \lambda = \frac{k}{r} - r $. We remark that only the $ n $-dimensional cylinder 
 $S^k (\sqrt k ) \times \mathbb{R}^{n-k}$, for $ 1 \leq k \leq n-1$,  is the self-shrinker of mean curvature flow. 
\end{Ex}

Let $ X(t) : M^n \times ( -\ep , \ep ) \rightarrow \mathbb{R}^{n+1} $ be a variation of $ X : M^n \rightarrow \mathbb{R}^{n+1} $. 
The weighted volume $ V : (-\ep , \ep ) \rightarrow \mathbb{R} $ is defined in \cite{CW} as follows:
\begin{eqnarray*}
V(t) = \int_M \innprod{ X(t)}{N} e^{-\frac{|X|}{2}} \dm .
\end{eqnarray*}
In \cite{CW}, Cheng and Wei have proved  that $ X : M^n \rightarrow \mathbb{R}^{n+1} $ is a critical point of the weighted area functional $ A(t) $ for the weighted volume-preserving variations if and only if $ X : M^n \rightarrow \mathbb{R}^{n+1} $ is a $ \lambda $-hypersurface. 
For further properties of $ \lambda $-hypersurfaces in details, see \cite{CW}. 
 
 In \cite{COW}, Cheng, Ogata and Wei have studied the rigidity theorem of complete $ \lambda $-hypersurfaces with a pointwise pinching
 condition. The following theorem is proved:
 \begin{Thm}
Let $X:M^n \rightarrow \mathbb{R}^{n+1}$ be an $n$-dimensional complete proper $\lambda$-hypersurface
in the Euclidean space $\mathbb{R}^{n+1}$.
If the squared norm $S$ of the second fundamental form and the mean curvature $H$ of $X:M^n \rightarrow \mathbb{R}^{n+1}$ satisfies
\begin{equation}\label{eq:001}
\biggl(\sqrt {S-\dfrac{H^2}n}+|\lambda| \dfrac{n-2}{2\sqrt{n(n-1)}}\biggl)^2+ \frac1n(H-\lambda)^2\leq 1+\dfrac{n\lambda^2}{4(n-1)},
\end{equation}
then $X:M^n \rightarrow \mathbb{R}^{n+1}$ is isometric to one of the following:
\begin{enumerate}
\item the sphere $S^n(r)$ with radius $0<r\leq\sqrt{n}$,
\item the Euclidean space  $\mathbb{R}^{n}$,
\item the cylinder $S^1(r)\times \mathbb{R}^{n-1}$ with radius $r>0$ and $n=2$ or with radius $r\geq 1$ and $n>2$,
\item the cylinder $S^{n-1}(r)\times \mathbb{R}$ with radius $r>0$ and $n=2$ or with radius $r\leq\sqrt{n-1}$ and $n>2$,
\item the cylinder $S^{k}(\sqrt{k})\times   \mathbb{R}^{n-k}$ for $2\leq k\leq n-2$.
\end{enumerate}
 \end{Thm}
  
In this paper, we study a global pinching theorem of complete $\lambda$-hypersurfaces in $ \mathbb{R}^{n+1} $. 
We prove  the following: 
\begin{Thm} \label{main}
Let $ X : M^n \rightarrow \mathbb{R}^{n+1} $ be an $n$-dimensional complete proper $\lambda $-hypersurface in $ \mathbb{R}^{n+1} $
with $ n \geq 3 $.  
If $ X : M^n \rightarrow \mathbb{R}^{n+1} $ satisfies 
\begin{equation}
\begin{aligned}
&\left( \int_M \left| \frac{n |\lambda| (n-2) }{ 2 \sqrt{n(n-1)} } B^\frac{1}{2} 
+ \frac{n}{2} B + \frac{n^2 -2n + 2}{2n(n-1)} H^2 - \frac{n+2}{2n} \lambda H \right|^\frac{n}{2} \dm \right)^\frac{2}{n}  \\
&< \frac{n-2}{n } k(n)^{-1} , 
\end{aligned}
\end{equation}
then $ X : M^n \rightarrow \mathbb{R}^{n+1} $ is isometric to the Euclidean space $ \mathbb{R}^n $ or the sphere $ S^n (r) $ with 
\begin{equation}
\left| r^2 - \frac{(3n-4)n}{(n-1)(n+2)} \right| 
< \frac{(n-2)^3}{4^{2(n+1)} n^{\frac{2}{n} +1} (n-1) (3n-4) (n+2) } \left( \frac{\omega_{n-1} }{ \omega_n } \right)^\frac{2}{n},
\end{equation} 
where  $ B = S - \dfrac{H^2}{n} $, 
$k(n) = \frac{2 \cdot 4^{2(n+1)} (n-1) (3n-4) }{ (n-2)^2 } \left( \frac{n}{\omega_{n-1} } \right)^\frac{2}{n}$, and 
  $\omega_{k} $ denotes the area  of the $k $-dimensional unit sphere $S^{k}(1) $. 
\end{Thm}
\vskip2mm
\noindent
{\bf Acknowledgement}. I would like to thank my adviser professor Qing-Ming Cheng for continuous guidance and support.

\section{the Sobolev inequality}
In order to prove our theorem, the following  Sobolev inequality in \cite{MS} plays a very important rule.

\begin{Thm} \label{Sobolev}
Let $ X : M^n \rightarrow \mathbb{R}^{n+1} $ be an $n$-dimensional hypersurface in $ \mathbb{R}^{n+1} $. For  any Lipschitz function $f\geq 0$ 
with compact support on $M $, 
\begin{eqnarray} \label{sobolev}
\left( \int_M f^\frac{n}{n-1} \dm \right)^\frac{n-1}{n} \leq C_n \int_M \left\{ |\nabla f| + | H | f \right\} \dm
\end{eqnarray}
holds,
where  
\begin{eqnarray*}
C_n = 4^{n+1} \left( \frac{n}{\omega_{n-1} } \right)^\frac{1}{n} .
\end{eqnarray*}
\end{Thm}
From the above theorem, we have the following corollary: 
\vskip2mm
\noindent
{\bf Corollary 2.1.}{ \it 
Let $ X : M^n \rightarrow \mathbb{R}^{n+1} $ be an $n$-dimensional hypersurface in $ \mathbb{R}^{n+1} $. For  any Lipschitz function 
$f\geq 0$ with compact support on $M $, 
\begin{eqnarray} \label{sob}
k(n)^{-1} \left( \int_M f^\frac{2n}{n-2} \dm \right)^\frac{n-2}{n} \leq \int_M |\nabla f|^2 \dm + \frac{n-2}{2(n-1)} \int_M H^2 f^2 \dm 
\end{eqnarray}
holds,
where  
\begin{eqnarray*}
k(n) = \frac{2 \cdot 4^{2(n+1)} (n-1) (3n-4) }{ (n-2)^2 } \left( \frac{n}{\omega_{n-1} } \right)^\frac{2}{n} .
\end{eqnarray*}
}

\paragraph{Proof}
Replacing $ f$ in  the theorem \ref{Sobolev} with  $f^\frac{2(n-1)}{n-2} $,  we get  
\begin{eqnarray*}
\left( \int_M f^\frac{2n}{n-2} \dm \right)^\frac{n-1}{n}
&\leq& C_n \int_M \left\{ |\nabla f^\frac{2(n-1)}{n-2} | + |H| f^\frac{2(n-1)}{n-2} \right\} \dm \\
&=& C_n \int_M \left\{ \frac{2(n-1)}{n-2} f^\frac{n}{n-2} |\nabla f| + |H| f^\frac{2(n-1)}{n-2} \right\} \dm .
\end{eqnarray*}
By H\"{o}lder's inequality, we have
\begin{eqnarray*}
\int_M f^\frac{n}{n-2} |\nabla f| \dm 
&\leq& \left( \int_M f^\frac{2n}{n-2} \dm \right)^\frac{1}{2} \left( \int_M |\nabla f|^2 \dm \right)^\frac{1}{2} , \\
\int_M |H| f^\frac{2(n-1)}{n-2} \dm &\leq& \left( \int_M f^\frac{2n}{n-2} \dm \right)^\frac{1}{2} \left( \int_M H^2 f^2 \dm \right)^\frac{1}{2} .
\end{eqnarray*}
Hence,  
\begin{eqnarray*}
\left( \int_M f^\frac{2n}{n-2} \dm \right)^\frac{n-1}{n} 
&\leq& C_n \frac{2(n-1)}{n-2} \left( \int_M f^\frac{2n}{n-2} \dm \right)^\frac{1}{2} \left( \int_M |\nabla f|^2 \dm \right)^\frac{1}{2} \\
&~& + C_n \left( \int_M f^\frac{2n}{n-2} \dm \right)^\frac{1}{2} \left( \int_M H^2 f^2 \dm \right)^\frac{1}{2} .
\end{eqnarray*}
Therefore 
\begin{eqnarray*}
\left( \int_M f^\frac{2n}{n-2} \dm \right)^\frac{n-2}{2n}
&\leq& C_n \Bigg\{ \frac{2(n-1)}{n-2} \left( \int_M |\nabla f|^2 \dm \right)^\frac{1}{2}  + \left( \int_M H^2 f^2 \dm \right)^\frac{1}{2} \Bigg\} .
\end{eqnarray*}
According to   $ \left( \int_M f^\frac{2n}{n-2} \dm \right)^\frac{n-2}{2n} \geq 0 $,  we have
\begin{eqnarray*}
 \left( \int_M f^\frac{2n}{n-2} \dm \right)^\frac{n-2}{n} 
&\leq& C_n^{~2} \Bigg\{ \frac{4(n-1)^2}{(n-2)^2}  \int_M |\nabla f|^2 \dm +  \int_M H^2 f^2 \dm  \\
&~&  + \frac{4(n-1)}{n-2} \left( \int_M |\nabla f|^2 \dm \right)^\frac{1}{2} \left( \int_M H^2 f^2 \dm \right)^\frac{1}{2} \Bigg\} .
\end{eqnarray*}
Because of  $ \int_M |\nabla f|^2 \dm \geq 0 $ and $ \int_M H^2 f^2 \dm \geq 0 $,  we get
\begin{eqnarray*}
&~& \left( \int_M f^\frac{2n}{n-2} \dm \right)^\frac{n-2}{n} \\
&\leq& C_n^{~2} \Bigg\{ \frac{4(n-1)^2}{(n-2)^2}  \int_M |\nabla f|^{2} \dm +  \int_M H^2 f^2 \dm  \\
&~&  + \frac{2(n-1)}{n-2} \left( \int_M |\nabla f|^2 \dm +  \int_M H^2 f^2 \dm \right) \Bigg\} \\
&=& C_n^{~2} \frac{2(n-1)}{n-2} \left(  \frac{2(n-1)}{n-2} +1\right) \left\{ \int_M |\nabla f|^2 \dm + \frac{n-2}{2(n-1)} \int_M H^2 f^2 \dm \right\}.
\end{eqnarray*}
Let $ k(n) = C_n^{~2} \frac{2(n-1)}{n-2} \left(  \frac{2(n-1)}{n-2} +1\right) $, then we get 
\begin{eqnarray*} 
k(n)^{-1} \left( \int_M f^\frac{2n}{n-2} \dm \right)^\frac{n-2}{n} \leq \int_M |\nabla f|^2 \dm + \frac{n-2}{2(n-1)} \int_M H^2 f^2 \dm .
\end{eqnarray*}
\fin

\section{Proof of our  global pinching theorem}

\vskip2mm
\noindent
In order to prove  the theorem \ref{main}, we prepare  several  lemmas. 
For the differential operator $\cL$ defined by 
$$
\cL f=\Delta f-\langle \nabla f, X\rangle = div ( e^{-\frac{|X|^2}{2}} \nabla f) e^\frac{|X|^2}{2}, 
$$
where $\Delta$ and $\nabla$ denote the Laplace operator and the gradient operator. 
In \cite{COW}, we have proved the following lemma. 

\begin{Lem} \label{L4}
For $ B= S - \frac{H^2}{n} $, we have  
\begin{eqnarray} \label{l4}
\frac{1}{2} \cL B \geq - \frac{|\lambda| (n-2) }{ \sqrt{n(n-1)} } B^\frac{3}{2} + B - B^2 - \frac{1}{n} H^2 B + \frac{2 \lambda }{n} HB .
\end{eqnarray}
\end{Lem}
Define a function $\rho$ by 
$$
\rho=e^{-\frac{|X|^2}2}.
$$
\begin{Lem} \label{L5}
For  any  smooth function $\eta$ with compact support on $ M $ and an arbitrary positive constant $ \ep $,  we have 
\begin{eqnarray} \label{l5}
&~& \int_M \Bigg\{ \frac{|\lambda|(n-2) }{ \sqrt{n(n-1)} } B^{n+ \frac{1}{2}} \eta^2 \rho- B^n \eta^2 \rho + B^{n+1} \eta^2 \rho \\
&~&+ \frac{1}{n} H^2 B^n \eta^2 \rho - \frac{2\lambda}{n} H B^n \eta^2 \rho + \frac{1}{2 \ep} B^n |\nabla \eta|^2 \rho \Bigg\} \dm \nonumber \\
&\geq& \frac{n-1-\ep}{2} \int_M B^{n-2} \eta^2 |\nabla B|^2 \rho ~\dm . \nonumber
\end{eqnarray}
\end{Lem}

\paragraph{Proof}
Multiplying $ B^{n-1} \eta^2 \rho $ on both sides of (\ref{l4}) and taking integral,  we obtain
\begin{eqnarray*}
 0 
&\geq& \int_M \Bigg\{ - \frac{|\lambda| (n-2) }{ \sqrt{n(n-1)} }  B^{n + \frac{1}{2}} \eta^2 \rho + B^n \eta^2 \rho - B^{n+1} \eta^2 \rho \\
&~& - \frac{1}{n} H^2 B^n \eta^2 \rho + \frac{2 \lambda }{n} HB^n \eta^2 \rho - \frac{1}{2} \cL B \cdot B^{n-1} \eta^2 \rho \Bigg\} \dm .
\end{eqnarray*}
Since $\eta$ has compact support on $M$, according to Stokes theorem, we get
\begin{eqnarray*}
&~& -\frac{1}{2} \int_M \cL B \cdot B^{n-1} \eta^2 \rho ~\dm \\
&=& -\frac{1}{2} \int_M div ( \rho \nabla B) \cdot B^{n-1} \eta^2 \dm \\
&=& \frac{1}{2} \int_M \innprod{ \rho \nabla B}{ \nabla ( B^{n-1} \eta^2 ) } \dm \\ 
&=& \frac{n-1}{2} \int_M B^{n-2} \eta^2 | \nabla B|^2 \rho ~\dm + \int_M B^{n-1} \eta \innprod{ \nabla B}{ \nabla \eta } \rho ~\dm.
\end{eqnarray*}
Moreover, for an arbitrary constant $ \ep > 0 $,   we have
\begin{eqnarray*}
\int_M B^{n-1} \eta \innprod{ \nabla B}{ \nabla \eta } \rho ~\dm
\geq -\frac{\ep}{2} \int_M B^{n-2} \eta^2 |\nabla B |^2 \rho ~\dm - \frac{1}{2 \ep}  \int_M B^n | \nabla \eta|^2 \rho ~\dm .
\end{eqnarray*}
Hence, we obtain
\begin{eqnarray*}
&~& \int_M \Bigg\{ \frac{|\lambda|(n-2) }{ \sqrt{n(n-1)} } B^{n+ \frac{1}{2}} \eta^2 \rho- B^n \eta^2 \rho + B^{n+1} \eta^2 \rho \\
&~& + \frac{1}{n} H^2 B^n \eta^2 \rho - \frac{2\lambda}{n} H B^n \eta^2 \rho + \frac{1}{2 \ep} B^n |\nabla \eta|^2 \rho \Bigg\} \dm \nonumber \\
&\geq& \frac{n-1-\ep}{2} \int_M B^{n-2} \eta^2 |\nabla B|^2 \rho ~\dm .
\end{eqnarray*}
\fin

\begin{Lem} \label{L6}
Putting $f := B^\frac{n}{2} \eta \rho^\frac{1}{2} $, we know that 
\begin{eqnarray} \label{l61}
&~& \int_M |\nabla f|^2 \dm  \\
&=& \frac{n^2}{4} \int_M B^{n-2} \eta^2 |\nabla B|^2 \rho ~\dm + \int_M B^n |\nabla \eta|^2 \rho ~\dm + n \int_M B^{n-1} \eta \innprod{ \nabla B}{ \nabla \eta } \rho ~\dm \nonumber \\
&~& - \frac{1}{4} \int_M | X^\top |^2 B^n \eta^2 \rho ~\dm + \frac{\lambda}{2} \int_M \innprod{ X}{ N} B^n \eta^2 \rho ~\dm - \frac{1}{2} |X^\bot|^2 B^n \eta^2 \rho ~\dm \nonumber \\
&~& + \frac{n}{2} \int_M B^n \eta^2 \rho ~\dm  \nonumber
\end{eqnarray}
and 
\begin{eqnarray} \label{l62}
&~& \frac{1}{2} \int_M H^2 f^2 \dm  \\
&=& \frac{ \lambda^2}{2} \int_M B^n \eta^2 \rho ~\dm - \lambda \int_M \innprod{ X}{N} B^n \eta^2 \rho ~\dm + \frac{1}{2} \int_M |X^\bot|^2 B^n \eta^2 \rho ~\dm . \nonumber
\end{eqnarray} 
hold.
\end{Lem}

\paragraph{Proof}
Calculating the left hand side of (\ref{l61}), we know
\begin{eqnarray*}
\int_M | \nabla f|^2 \dm 
&=& \int_M | \nabla ( B^\frac{n}{2} \eta ) |^2 \rho ~\dm + \int_M B^n \eta^2 | \nabla \rho^\frac{1}{2} |^2 \dm \\
&~& + 2 \int_M B^\frac{n}{2} \eta \rho^\frac{1}{2} \innprod{ \nabla ( B^\frac{n}{2} \eta ) }{ \nabla \rho^\frac{1}{2} } \dm .
\end{eqnarray*}
Putting
\begin{eqnarray*}
T_1 &:=& \int_M | \nabla ( B^\frac{n}{2} \eta ) |^2 \rho ~\dm , \\
T_2 &:=& \int_M B^n \eta^2 | \nabla \rho^\frac{1}{2} |^2 \dm , \\
T_3 &:=& 2 \int_M B^\frac{n}{2} \eta \rho^\frac{1}{2} \innprod{ \nabla ( B^\frac{n}{2} \eta ) }{ \nabla \rho^\frac{1}{2} } \dm .
\end{eqnarray*}
\begin{eqnarray*}
T_1 
&=& \frac{n^2}{4} \int_M B^{n-2} \eta^2 |\nabla B|^2 \rho ~\dm + \int_M B^n |\nabla \eta|^2 \rho ~\dm \\
&~& + n \int_M B^{n-1} \eta \innprod{ \nabla B}{ \nabla \eta } \rho ~\dm .
\end{eqnarray*}
Because of  $ | \nabla \rho^\frac{1}{2} |^2 = \frac{1}{4} | X^\top|^2 \rho $ and $ \Delta X = H N $, we have 
\begin{eqnarray*}
T_2 = \frac{1}{4} \int_M |X^\top|^2 B^n \eta^2 \rho ~\dm 
\end{eqnarray*}
and
\begin{eqnarray*}
\Delta \rho 
&=& |X^\top|^2\rho - \innprod{ \Delta X }{ X } \rho - n \rho \\
&=& | X^\top|^2 \rho - \lambda \innprod{ X}{ N } \rho + | X^\bot|^2 \rho - n \rho .
\end{eqnarray*}
Hence, 
\begin{eqnarray*}
 T_3 
&=& \frac{1}{2} \int_M \innprod{ \nabla ( B^n \eta^2 ) }{ \nabla \rho } \dm\\
&=& - \frac{1}{2} \int_M B^n \eta^2 \cdot \Delta \rho ~\dm \\
&=& - \frac{1}{2} \int_M |X^\top|^2 B^n \eta^2 \rho ~\dm + \frac{ \lambda }{ 2} \int_M \innprod{ X}{N} B^n \eta^2 \rho ~\dm  \\
&~& - \frac{1}{2}\int_M  |X^\bot|^2 B^n \eta^2 \rho ~\dm+ \frac{n}{2} \int_M B^n \eta^2 \rho ~\dm .
\end{eqnarray*}
Therefore, we get 
\begin{eqnarray*}
&~& \int_M |\nabla f|^2 \dm  \\
&=& \frac{n^2}{4} \int_M B^{n-2} \eta^2 |\nabla B|^2 \rho ~\dm + \int_M B^n |\nabla \eta|^2 \rho ~\dm + n \int_M B^{n-1} \eta \innprod{ \nabla B}{ \nabla \eta } \rho ~\dm \nonumber \\
&~& - \frac{1}{4} \int_M | X^\top |^2 B^n \eta^2 \rho ~\dm + \frac{\lambda}{2} \int_M \innprod{ X}{ N} B^n \eta^2 \rho ~\dm - \frac{1}{2} |X^\bot|^2 B^n \eta^2 \rho ~\dm \\
&~& + \frac{n}{2} \int_M B^n \eta^2 \rho ~\dm.
\end{eqnarray*}
From $ H = \lambda - \innprod{ X}{ N} $,  we get
\begin{eqnarray*}
\frac{1}{2} \int_M H^2 f^2 \dm 
&=& \frac{1}{2} \int_M ( \lambda - \innprod{ X}{N} )^2 B^n \eta^2 \rho ~\dm \\
&=& \frac{\lambda^2 }{2} \int_M B^n \eta^2 \rho ~\dm - \lambda \int_M \innprod{X}{N} B^n \eta^2 \rho ~\dm \\
&~& + \frac{1}{2} \int_M |X^\bot|^2 B^n \eta^2 \rho ~\dm .
\end{eqnarray*}
\fin

\begin{Lem} \label{L7}
For an arbitrary constant $ \delta > 0$,  we have
\begin{eqnarray} \label{l7}
&~& k(n)^{-1} \left( \int_M f^\frac{2n}{n-2} \dm \right)^\frac{n-2}{n}  \\
&\leq& \frac{(1+\delta)n^2}{4} \int_M B^{n-2} \eta^2 | \nabla B|^2 \rho ~\dm + \left( 1+ \frac{1}{\delta} \right) \int_M B^n | \nabla \eta|^2 \rho ~\dm  \nonumber \\
&~& - \frac{1}{2(n-1) } \int_M H^2 B^n \eta^2 \rho ~\dm +\frac{\lambda}{2} \int_M H B^n \eta^2 \rho ~\dm + \frac{n}{2} \int_M B^n \eta^2 \rho ~\dm , \nonumber
\end{eqnarray}
where $ k(n)$ is the assertion of the Corollary 2.1. 
\end{Lem}

\paragraph{Proof}
From Corollary 2.1, we  have, for any function $f$ with compact support on $M$, 
\begin{eqnarray*}
k(n)^{-1} \left( \int_M f^\frac{2n}{n-2} \dm \right)^\frac{n-2}{n} 
&\leq& \int_M |\nabla f|^2 \dm + \frac{n-2}{2(n-1)} \int_M H^2 f^2 \dm \\
&=& \int_M |\nabla f|^2 \dm + \frac{1}{2} \int_M H^2 f^2 \dm - \frac{1}{2(n-1)} \int_M H^2 f^2 \dm . 
\end{eqnarray*}
Taking  $ f=B^\frac{n}{2} \eta \rho^\frac{1}{2} $, 
from Lemma \ref{L6},  we infer
\begin{eqnarray*}
&~& k(n)^{-1} \left( \int_M f^\frac{2n}{n-2} \dm \right)^\frac{n-2}{n} \\
&\leq& \frac{n^2}{4} \int_M B^{n-2} \eta^2 |\nabla B|^2 \rho ~\dm + \int_M B^n |\nabla \eta|^2 \rho ~\dm + n \int_M B^{n-1} \eta \innprod{ \nabla B}{ \nabla \eta } \rho ~\dm  \\
&~& - \frac{1}{4} \int_M | X^\top |^2 B^n \eta^2 \rho ~\dm + \frac{\lambda}{2} \int_M \innprod{ X}{ N} B^n \eta^2 \rho ~\dm - \frac{1}{2} |X^\bot|^2 B^n \eta^2 \rho ~\dm \\
&~& + \frac{n}{2} \int_M B^n \eta^2 \rho ~\dm +\frac{\lambda^2 }{2} \int_M B^n \eta^2 \rho ~\dm - \lambda \int_M \innprod{X}{N} B^n \eta^2 \rho ~\dm  \\
&~& + \frac{1}{2} \int_M |X^\bot|^2 B^n \eta^2 \rho ~\dm - \frac{1}{2(n-1)} \int_M H^2 B^n \eta^2 \rho ~ \dm \\
&\leq& \frac{n^2}{4} \int_M B^{n-2} \eta^2 |\nabla B|^2 \rho ~\dm + \int_M B^n |\nabla \eta|^2 \rho ~\dm + n \int_M B^{n-1} \eta \innprod{ \nabla B}{ \nabla \eta } \rho ~\dm  \\
&~& - \frac{\lambda}{2} \int_M \innprod{ X}{ N} B^n \eta^2 \rho ~\dm - \frac{1}{2(n-1)} \int_M H^2 B^n \eta^2 \rho ~\dm + \left( \frac{n}{2} + \frac{\lambda^2}{2} \right) \int_M B^n \eta^2 \rho ~\dm \\
&=& \frac{n^2}{4} \int_M B^{n-2} \eta^2 |\nabla B|^2 \rho ~\dm + \int_M B^n |\nabla \eta|^2 \rho ~\dm + n \int_M B^{n-1} \eta \innprod{ \nabla B}{ \nabla \eta } \rho ~\dm  \\
&~& + \frac{\lambda}{2} \int_M H B^n \eta^2 \rho ~\dm - \frac{1}{2(n-1)} \int_M H^2 B^n \eta^2 \rho ~\dm +  \frac{n}{2}  \int_M B^n \eta^2 \rho ~\dm .
\end{eqnarray*}
For an arbitrary constant $ \delta > 0 $,  we have
\begin{eqnarray*}
n \int_M B^{n-1} \eta \innprod{ \nabla B}{ \nabla \eta } \rho ~\dm
\leq \frac{\delta n^2}{4} \int_M B^{n-2} \eta^2 |\nabla B|^2 \rho ~\dm + \frac{1}{\delta} \int_M B^n |\nabla \eta|^2 \rho ~\dm .
\end{eqnarray*}
Hence,  we get
\begin{eqnarray*}
&~& k(n)^{-1} \left( \int_M f^\frac{2n}{n-2} \dm \right)^\frac{n-2}{n} \nonumber \\
&\leq& \frac{(1+\delta)n^2}{4} \int_M B^{n-2} \eta^2 | \nabla B|^2 \rho ~\dm + \left( 1+ \frac{1}{\delta} \right) \int_M B^n | \nabla \eta|^2 \rho ~\dm  \nonumber \\
&~&- \frac{1}{2(n-1) } \int_M H^2 B^n \eta^2 \rho ~\dm +\frac{\lambda}{2} \int_M H B^n \eta^2 \rho ~\dm + \frac{n}{2} \int_M B^n \eta^2 \rho ~\dm .
\end{eqnarray*}
\fin

\paragraph{{\it Proof of Theorem} \ref{main}}
If  $ B \not\equiv 0 $ holds, we can choose $\eta$ such that, for $ f = B^{\frac12} \eta \rho^{\frac12} $, 
$$
\left( \int_M f^\frac{2n}{n-2} \dm \right)^\frac{n-2}{n}\neq 0.
$$
From Lemma \ref{L5} and Lemma \ref{L7}, then for arbitrary constants $ \ep > 0$ and $ \delta > 0 $, 
\begin{eqnarray*}
&~& k(n)^{-1} \left( \int_M f^\frac{2n}{n-2} \dm \right)^\frac{n-2}{n} \\
&\leq& \frac{(1+\delta)n^2}{2} \frac{1}{n-1-\ep} 
\int_M \Bigg\{ \frac{|\lambda|(n-2) }{ \sqrt{n(n-1)} } B^{n+ \frac{1}{2}} \eta^2 \rho- B^n \eta^2 \rho 
+ B^{n+1} \eta^2 \rho \\
&~& + \frac{1}{n} H^2 B^n \eta^2 \rho - \frac{2\lambda}{n} H B^n \eta^2 \rho 
+ \frac{1}{2 \ep} B^n |\nabla \eta|^2 \rho \Bigg\} \dm \\
&~& + \left( 1+ \frac{1}{\delta} \right) \int_M B^n | \nabla \eta|^2 \rho ~\dm 
- \frac{1}{2(n-1) } \int_M H^2 B^n \eta^2 \rho ~\dm \nonumber \\
&~& +\frac{\lambda}{2} \int_M H B^n \eta^2 \rho ~\dm + \frac{n}{2} \int_M B^n \eta^2 \rho ~\dm .
\end{eqnarray*}
Letting  $1 + \delta = \dfrac{n-1+\ep}{n} $, then, we derive 
\begin{eqnarray*}
&~& k(n)^{-1} \left( \int_M f^\frac{2n}{n-2} \dm \right)^\frac{n-2}{n} \\
&\leq& \frac{n-1+\ep}{n-1-\ep} \int_M \Bigg\{ \frac{n|\lambda|(n-2)}{2\sqrt{n(n-1)} } B^{n+\frac{1}{2} } \eta^2 \rho
 + \frac{n}{2} B^{n+1} \eta^2 \rho   \\
&~& + \frac{1}{2} \left(1 - \frac{1}{(n-1)} \frac{n-1-\ep}{n-1+\ep} \right) H^2 B^n \eta^2 \rho
+ \left( - 1 + \frac{1}{2} \frac{n-1-\ep}{n-1+\ep} \right) \lambda H B^n \eta^2 \rho \Bigg\} \dm \\
&~& +\frac{n}{2} \left( - \frac{n-1+\ep}{n-1-\ep} + 1 \right) \int_M B^n \eta^2 \rho ~\dm 
+ C(n,\ep) \int_M B^n |\nabla \eta|^2 \rho ~\dm \\
&\leq& \frac{n-1+\ep}{n-1-\ep} \int_M \Bigg\{ \frac{n|\lambda|(n-2)}{2\sqrt{n(n-1)} } B^{n+\frac{1}{2} } \eta^2 \rho 
+ \frac{n}{2} B^{n+1} \eta^2 \rho   \\
&~& + \frac{1}{2} \left(1 - \frac{1}{(n-1)} \frac{n-1-\ep}{n-1+\ep} \right) H^2 B^n \eta^2 \rho
+ \left( - 1 + \frac{1}{2} \frac{n-1-\ep}{n-1+\ep} \right) \lambda H B^n \eta^2 \rho \Bigg\} \dm \\
&~& + C(n,\ep) \int_M B^n |\nabla \eta|^2 \rho ~\dm, 
\end{eqnarray*}
where $ C(n,\ep) $ is a positive constant only depending on $ n $ and $ \ep $. 
From $ f^2 = B^n \eta^2 \rho $ and  using H\"{o}lder's inequality,  we obtain
\begin{eqnarray*}
&~& k(n)^{-1} \left( \int_M f^\frac{2n}{n-2} \dm \right)^\frac{n-2}{n} \\
&\leq& \frac{n-1+\ep}{n-1-\ep} \Bigg( \int_M \bigg| \frac{n|\lambda|(n-2)}{2\sqrt{n(n-1)} } B^\frac{1}{2}   + \frac{n}{2} B   + \frac{1}{2} \left(1 - \frac{1}{(n-1)} \frac{n-1-\ep}{n-1+\ep} \right) H^2 \\
&~& + \left( - 1 + \frac{1}{2} \frac{n-1-\ep}{n-1+\ep} \right) \lambda H  \bigg|^\frac{n}{2} \dm \Bigg)^\frac{2}{n}  \left( \int_M f^\frac{2n}{n-2} \dm \right)^\frac{n-2}{n} \\
&~&+C(n,\ep) \int_M B^n |\nabla \eta|^2 \rho ~\dm.
\end{eqnarray*}
Therefore, we have

\begin{eqnarray*}
&~& k(n)^{-1}  \\
&\leq& \frac{n-1+\ep}{n-1-\ep} \Bigg( \int_M \bigg| \frac{n|\lambda|(n-2)}{2\sqrt{n(n-1)} } B^\frac{1}{2}   + \frac{n}{2} B   + \frac{1}{2} \left(1 - \frac{1}{(n-1)} \frac{n-1-\ep}{n-1+\ep} \right) H^2 \\
&~& + \left( - 1 + \frac{1}{2} \frac{n-1-\ep}{n-1+\ep} \right) \lambda H  \bigg|^\frac{n}{2} \dm \Bigg)^\frac{2}{n}   \\
&~&+C(n,\ep) \dfrac{\int_M B^n |\nabla \eta|^2 \rho ~\dm}{\left( \int_M f^\frac{2n}{n-2} \dm \right)^\frac{n-2}{n}}.
\end{eqnarray*}
Since $X: M^n \rightarrow \mathbb{R}^{n+1} $ is proper, it is proved by Cheng and Wei in \cite{CW} that 
$X: M^n \rightarrow \mathbb{R}^{n+1}$ has at most polynomial area growth. Hence, we know that 
$$
\int_M B^n  \rho ~\dm<\infty.
$$
Taking $\eta=\phi(\frac{|X|}{r})$ for any $r>0$, where $\phi$ is a nonnegative function on $[0, \infty)$ such that 
\begin{equation*}
\phi(t)=\begin{cases}
1, &\text{ if  $t\in [0,1]$} \\
0, &\text{ if $t\in [2, \infty)$}
\end{cases}
\end{equation*}
and $|\phi^{\prime}|\leq c$ for some absolute constant. Taking $r\to\infty$, we have
$$
\int_M B^n |\nabla \eta|^2 \rho ~\dm\to 0.
$$ 
Therefore, we get 
\begin{eqnarray*}
&~& k(n)^{-1}  \\
&\leq& \frac{n-1+\ep}{n-1-\ep} \Bigg( \int_M \bigg| \frac{n|\lambda|(n-2)}{2\sqrt{n(n-1)} } B^\frac{1}{2}   + \frac{n}{2} B   + \frac{1}{2} \left(1 - \frac{1}{(n-1)} \frac{n-1-\ep}{n-1+\ep} \right) H^2 \\
&~& + \left( - 1 + \frac{1}{2} \frac{n-1-\ep}{n-1+\ep} \right) \lambda H  \bigg|^\frac{n}{2} \dm \Bigg)^\frac{2}{n}.
\end{eqnarray*}
Letting $ \ep \rightarrow 1$, we obtain
\begin{equation*}
\begin{aligned}
& k(n)^{-1}  \\
&\leq \frac{n}{n-2} \biggl( \int_M \bigg| \frac{n|\lambda|(n-2)}{2\sqrt{n(n-1)} } B^\frac{1}{2}  + \frac{n}{2} B 
+ \frac{n^2-2n+2}{2n(n-1)}  H^2  -  \frac{n+2}{2n} \lambda H  \bigg|^\frac{n}{2} \dm \biggl)^\frac{2}{n}\\
&< \frac{n}{n-2} \cdot \frac{n-2}{n  } k(n)^{-1}= k(n)^{-1} .
\end{aligned}
\end{equation*}
It is a  contradiction. 
Thus, we have $B = S-\frac{H^2}{n} \equiv 0$, that is,  $X: M^n \rightarrow \mathbb{R}^{n+1} $ is totally umbilical. Hence, we know  that 
 $X: M^n \rightarrow \mathbb{R}^{n+1} $ is isomeric to $ \mathbb{R}^n $ or  a sphere $ S^n(r)$ with radius $r$, which satisfies (1.3) from (1.2).
\fin

\vskip1cm

Shiho Ogata 

Department of Applied Mathematics, Graduate School of Sciences, 

Fukuoka  University, Fukuoka 814-0180,  Japan

\end{document}